\documentclass[a4paper,11pt]{article}
\usepackage{amsfonts,amsmath,graphicx,calc,enumerate,url}
\usepackage[sort&compress,numbers]{natbib}
\usepackage[lmargin=37mm,rmargin=37mm,tmargin=33mm,bmargin=33mm]{geometry}

\usepackage{amsthm}

\theoremstyle{plain}

\newtheorem{theorem}{Theorem}

\newtheorem{lemma}{Lemma}
\newtheorem{corollary}{Corollary}

\newtheorem{proposition}{Proposition}

\theoremstyle{definition}
\newtheorem{conjecture}{Conjecture}

\newtheorem{open}{Open Problem}

\newcommand{\thmlabel}[1]{\label{thm:#1}}
\newcommand{\thmref}[1]{Theorem~\ref{thm:#1}}

\newcommand{\lemlabel}[1]{\label{lem:#1}}
\newcommand{\lemref}[1]{Lemma~\ref{lem:#1}}

\newcommand{\eqnlabel}[1]{\label{eqn:#1}}
\newcommand{\eqnref}[1]{\eqref{eqn:#1}}

\newcommand{\figlabel}[1]{\label{fig:#1}}
\newcommand{\figref}[1]{Figure~\ref{fig:#1}}

\newcommand{\seclabel}[1]{\label{sec:#1}}
\newcommand{\secref}[1]{Section~\ref{sec:#1}}

\newcommand{\mySection}[2]{\section{\boldmath\textsc{#1}}\seclabel{#2}}
\renewcommand{\mySection}[2]{\section{#1}\seclabel{#2}}
\newcommand{\mySubSection}[2]{\subsection{\boldmath #1}\seclabel{#2}}

\newcommand{\corlabel}[1]{\label{cor:#1}}
\newcommand{\corref}[1]{Corollary~\ref{cor:#1}}

\newcommand{\proplabel}[1]{\label{prop:#1}}
\newcommand{\propref}[1]{Proposition~\ref{prop:#1}}

\newcommand{\openlabel}[1]{\label{open:#1}}
\newcommand{\openref}[1]{Open Problem~\ref{open:#1}}

\DeclareFontShape{OT1}{cmr}{b}{sc}
{<5><6><7><8><9><10><12><10.95><14.4><17.28><20.74><24.88>cmbcsc10}{}
\DeclareFontShape{OT1}{cmr}{bx}{sc}
{<5><6><7><8><9><10><12><10.95><14.4><17.28><20.74><24.88>cmbcsc10}{}

\DeclareFontShape{OT1}{cmr}{b}{tt}
{<5><6><7><8><9><10><12><10.95><14.4><17.28><20.74><24.88>cmbtt10}{}
\DeclareFontShape{OT1}{cmr}{bx}{tt}
{<5><6><7><8><9><10><12><10.95><14.4><17.28><20.74><24.88>cmbtt10}{}

\newcommand{\R}{\ensuremath{\mathbb{R}}}

\newcommand{\FancyFootnotes}{\renewcommand{\thefootnote}{\fnsymbol{footnote}}}
\newcommand{\FancyFootnotesOff}{\renewcommand{\thefootnote}{\arabic{footnote}}}

\newcommand{\e}{\ensuremath{\boldsymbol{e}}}

\newcommand{\paran}[1]{\textup{(}#1\textup{)}}

\newcommand{\SET}[1]{\ensuremath{\protect\left\{#1\right\}}}

\newcommand{\BRACKET}[1]{\ensuremath{\protect\left(#1\right)}}
\newcommand{\BFRAC}[2]{\ensuremath{\protect\left(\frac{#1}{#2}\right)}}
\newcommand{\bfrac}[2]{\ensuremath{\protect(\frac{#1}{#2})}}

\newcommand{\ceil}[1]{\ensuremath{\protect\lceil#1\rceil}}

\newcommand{\floor}[1]{\ensuremath{\protect\lfloor#1\rfloor}}

\newcommand{\half}{\ensuremath{\protect\tfrac{1}{2}}}

\newcommand{\etal}{~\emph{et~al.}~}

\newcommand{\Oh}[1]{\ensuremath{\protect\mathcal{O}(#1)}}

\newcommand{\Figure}[4][htb]{
\begin{figure}[#1]
	\vspace*{1ex}
	\begin{center}#3\end{center}
	\vspace*{-1ex}
	\caption{\figlabel{#2}#4}
\end{figure}
}

\newlength{\marginboxwidth}

\usepackage{pifont}

\newcommand{\mequiv}[3]{\ensuremath{#1\equiv #2}}
\newcommand{\GGG}{\ensuremath{\mathcal{G}}}

\newcommand{\sn}[1]{\ensuremath{\textup{\textsf{sn}}(#1)}}
\newcommand{\csn}[1]{\ensuremath{\textup{\textsf{csn}}(#1)}}
\newcommand{\bw}[1]{\ensuremath{\textup{\textsf{bw}}(#1)}}
\newcommand{\ppw}[1]{\ensuremath{\textup{\textsf{ppw}}(#1)}}

\begin{document}

\date{Submitted: July 28, 2005;\quad Revised: \today}


\title{\textbf{Graph Drawings with Few Slopes}\,\thanks{A preliminary version of this paper was published as: ``Really straight graph drawings.'' \emph{Proceedings of the 12th International Symposium on Graph Drawing} (GD '04), \emph{Lecture Notes in Computer Science} 3383:122--132, Springer, 2004.}\\[2ex] \emph{\large Dedicated to Godfried Toussaint on his 60th birthday.}\\[2ex]}

\FancyFootnotes

\author{
Vida Dujmovi{\'c} \footnotemark[2] \and 
Matthew Suderman \footnotemark[4] \and 
David R. Wood \footnotemark[5]
}

\footnotetext[2]{Department of Mathematics and Statistics, McGill University, Montr{\'e}al, Qu\'ebec, Canada (\texttt{vida@cs.mcgill.ca}). Supported by NSERC.}

\footnotetext[4]{McGill Centre for Bioinformatics, School of Computer Science, McGill University, Montr{\'e}al, Qu\'ebec, Canada (\texttt{msuder@cs.mcgill.ca}). Supported by NSERC.}

\footnotetext[5]{Departament de Matem{\`a}tica Aplicada II, Universitat Polit{\`e}cnica de Catalunya, Barcelona, Catalunya, Spain (\texttt{david.wood@upc.edu}). Supported by a Marie Curie 
Fellowship of the European Community under contract 023865, and by the 
projects MCYT-FEDER BFM2003-00368 and Gen.\ Cat 2001SGR00224.}

\maketitle

\begin{abstract}  
The \emph{slope-number} of a graph $G$ is the minimum number of distinct edge slopes in a straight-line drawing of $G$ in the plane. We prove that for $\Delta\geq5$ and all large $n$, there is a $\Delta$-regular $n$-vertex graph with slope-number at least $n^{1-\frac{8+\varepsilon}{\Delta+4}}$. This is the best known lower bound on the slope-number of a graph with bounded degree. We prove upper and lower bounds on the slope-number of complete bipartite graphs. We prove a general upper bound on the slope-number of an arbitrary graph in terms of its bandwidth. It follows that the slope-number of interval graphs, cocomparability graphs, and AT-free graphs is at most a function of the maximum degree. We prove that graphs of bounded degree and bounded treewidth have slope-number at most \Oh{\log n}. Finally we prove that every graph has a drawing with one bend per edge, in which the number of slopes is at most one more than the maximum degree. In a companion paper, planar drawings of graphs with few slopes are also considered. 
\end{abstract}

\newpage

\FancyFootnotesOff

\mySection{Introduction}{Introduction}

This paper studies straight-line drawings of graphs\footnote{We consider  undirected, finite, and simple graphs $G$ with  vertex set $V(G)$ and edge set $E(G)$. The  number of vertices and edges of $G$ are respectively denoted by $n=|V(G)|$ and $m=|E(G)|$. The minimum and maximum degrees of $G$ are respectively denoted by $\delta(G)$ and $\Delta(G)$.} in the plane with few distinct edge slopes\footnote{Consider a mapping of the vertices of a graph to distinct points in the plane. Now represent each edge by the closed line segment between its endpoints. Such a mapping is a (\emph{straight-line}) \emph{drawing} if the only vertices that each edge intersects is its own endpoints. The \emph{slope} of a line $L$ is the angle swept from the X-axis in an anticlockwise direction to $L$ (and is thus in $[0,\pi)$). The \emph{slope} of an edge or segment is the slope of the line that contains it. Of course two edges have the same \emph{slope} if and only if they are parallel. A \emph{crossing} in a drawing is a pair of edges that intersect at some point other than a common endpoint. A drawing is \emph{plane} if it has no crossings.}. \citet{WadeChu-CJ94} introduced this topic, and defined the \emph{slope-number} of a graph $G$ to be the minimum number of distinct edge slopes in a drawing of $G$. Let the \emph{convex slope-number} of $G$ be the minimum number of distinct edge slopes in a convex drawing\footnote{A drawing is \emph{convex} if all the vertices are on the convex hull, and no three vertices are collinear.} of $G$. Let \sn{G} and \csn{G} respectively be the slope-number and convex slope-number of $G$. By definition $\sn{G}\leq\csn{G}$ for every graph $G$. In this paper we prove lower and upper bounds on \sn{G} and \csn{G} for various (families of) graphs $G$. In a companion paper \citep{DESW-PlanarSlopesSegments}, planar drawings of graphs with few slopes are also considered. 

We start by considering some elementary lower bounds on the number of slopes. In a drawing of a graph, at most two edges incident to a vertex $v$ can have the same slope. Thus the edges incident to $v$ use at least $\half\deg(v)$ slopes. Hence the number of slopes is at least half the maximum degree. For some vertex $v$ on the convex hull of the drawing, every edge incident to $v$ has a distinct slope. Thus the number of slopes is at least the minimum degree. In a convex drawing, every edge incident to each vertex $v$ has a distinct slope. Thus the number of slopes is at least the maximum degree. Summarising:
\begin{equation}
\eqnlabel{LowerBoundSlopes} 
\text{(a) }\sn{G}\geq\half\Delta(G),
\quad
\text{(b) }\sn{G}\geq\delta(G),
\quad
\text{and}
\quad\text{(c) }\csn{G}\geq\Delta(G).
\end{equation}

Given these three lower bounds, it is natural ask whether there is a function $f$ such that $\sn{G}\leq f(\Delta(G))$ for every graph $G$. (A result by \citet{Malitz94a} implies that there is no such function $f$ for convex slope-number\footnote{The \emph{book thickness} of a graph $G$ is the minimum integer $k$ such that $G$ has a drawing in which each edge receives one of $k$ colours, and edges with the same colour do not cross; see \citep{DujWoo-DMTCS04}. Since parallel edges do not cross, the book thickness of $G$ is a lower bound on \csn{G}. \citet{Malitz94a} proved that there are $\Delta$-regular $n$-vertex graphs $G$ with book thickness $\Omega(\sqrt{\Delta}n^{1/2-1/\Delta})$. \citet{BMW-EJC06} proved the same result for all $\Delta\geq3$. Thus $\csn{G}\geq \Omega(\sqrt{\Delta}n^{1/2-1/\Delta})$.}.)\ This question was first posed in the conference version of this paper \citep{ReallyStraight-GD04}. It was subsequently solved in the negative for $\Delta\geq5$ independently by \citet{PachPal-EJC06} and \citet{BMW-EJC06} \footnote{The \emph{geometric thickness} of a graph $G$ is the minimum integer $k$ such that $G$ has a drawing in which each edge receives one of $k$ colours, and edges with the same colour do not cross; see \citep{DEH-JGAA00, DujWoo-GD05, Eppstein-AMS, DEK-SoCG04}. Since parallel edges do not cross, the geometric thickness of $G$ is a lower bound on \sn{G}. \citet{BMW-EJC06} proved that for all $\Delta\geq9$ and $\varepsilon>0$, for all sufficiently large $n>n(\Delta,\varepsilon)$, there exists a $\Delta$-regular $n$-vertex graph with geometric thickness at least 
$c\sqrt{\Delta}\,n^{1/2-4/\Delta-\varepsilon}$.}. The best bound, due to  \citet{PachPal-EJC06}, states that for all $\Delta\geq5$ and for all sufficiently large $n$, there exists an $n$-vertex graph $G$ with maximum degree $\Delta$ and slope-number 
\begin{equation*}
\sn{G}\;>\;n^{\frac{1}{2}-\frac{1}{\Delta-2}-\text{o}(1)}.
\end{equation*}
The first contribution of this paper is to prove an analogous lower bound of
\begin{equation*}
\sn{G}\;>\;n^{1-\frac{8+\varepsilon}{\Delta+4}}.
\end{equation*}
for all $\Delta\geq5$ (\secref{Degree}). This is the best known bound for all $\Delta\geq 9$. More importantly, our bound tends to $n$ for large $\Delta$, whereas the previous bounds by \citet{PachPal-EJC06} and \citet{BMW-EJC06} both tend to $\sqrt{n}$. 

The other main contributions of this paper establish graph families for which the slope-number is at most a function of the maximum degree. First, we consider the slope-number of complete $k$-partite graphs (\secref{Multipartite}). 

We then show that the slope-number is at most a function of the maximum degree for interval graphs\footnote{A graph $G$ is an \emph{interval graph} if one can assign to each vertex $v\in V(G)$ a closed interval $[L_v,R_v]\subset\mathbb{R}$ such that $vw\in E(G)$ if and only if $[L_v,R_v]\cap[L_w,R_w]\ne\emptyset$. The \emph{pathwidth} of a graph $G$ is the minimum $k$ such that $G$ is a spanning subgraph of an interval graph with no clique on $k+2$ vertices.}, cocomparability graphs\footnote{Let $\preceq$ be a partial order on a ground set $P$. The \emph{cocomparability} graph of $\preceq$ has vertex set $P$, where two vertices are adjacent if they are incomparable under $\preceq$. For example, every permutation graph is a cocomparability graph.}, and AT-free graphs\footnote{An \emph{asteroidal triple} in a graph is an independent set of three vertices such that each pair is joined by a path that avoids the neighborhood of the third. A graph is \emph{asteroidal triple-free} (or \emph{AT-free}) if it contains no asteroidal triple. AT-free graphs include interval, trapezoid, and cocomparability graphs.}. These results are established by first proving a general upper bound on the slope-number in terms of the bandwidth (\secref{Paths}). 

For graphs with bounded degree and bounded treewidth\footnote{A graph is \emph{chordal} if every induced cycle is a triangle. The \emph{treewidth} of a graph $G$ is the minimum integer $k$ such that $G$ is a subgraph of a chordal graph with no clique on $k+2$ vertices. This parameter is particularly important in algorithmic and structural graph theory; see \citep{Bodlaender-TCS98, Reed-AlgoTreeWidth03} for surveys. The treewidth of a graph is at most its pathwidth.}, we prove a \Oh{\log n} upper bound on the slope-number (\secref{Trees}). The proof is based on a result of independent interest: every tree $T$ has a drawing with $\Delta(T)-1$ slopes and $2k-1$ distinct edge lengths, where $k$ is the pathwidth of $T$. 

Our final contribution is to show that every graph $G$ has a drawing with $\Delta(G)+1$ slopes, if we allow one bend in each edge (\secref{Bends}).



\mySubSection{Related Research}{Related}

We now outline some related research from the literature. Drawings of lattices and posets with few slopes have been considered by \citet{FJ70}, Czyzowicz\etal\citep{Czy-JCTA91, CPRU-Order90, CPR-DM90} and \citet{Freese-04}.

\citet{ABH05} introduced the following slope parameter of graphs. Let $P\subset\R^2$ be a finite set of points in the plane. Let $S\subset\R\cup\{\infty\}$ be a set of slopes. Let $G(P,S)$ be the graph with vertex set $P$ where two points $v,w\in P$ are adjacent if and only if the slope of the line $\overline{vw}$ is in $S$. The \emph{slope parameter} of a graph $G$ is the minimum integer $k$ such that $G\cong G(P,S)$ for some point set $P$ and slope set $S$ with $|S|k$. This idea differs from our definition in that a clique can be represented by a set of collinear points. Amongst other results, \citet{ABH05} characterised the graphs with slope parameter $2$, and proved that the slope parameter of a tree $T$ equals $\Delta(T)$.

A famous result by  \citet{Ungar-JCTA82}, settling an open problem of \citet{Scott70}, states that $n$ non-collinear points determine at least $n-1$ distinct slopes. The configurations of $n$ points that determine exactly $n-1$ distinct slopes have been investigated by \citet{Jamison84, Jamison84a}. \citet{Jamison87} generalised the result of Ungar by proving that any set of non-collinear points has a spanning tree whose edges have distinct slopes.
\citet{Jamison87} conjectured that any set of points in general position has a spanning path whose edges have distinct slopes. In this direction, \citet{KP-DCG05} proved that every $n$-vertex caterpillar has a drawing on any $n$ prespecified points in general position such that no two edges have the same slope. 




Multi-dimensional graph drawings with few slopes are also of interest. Since an orthogonal projection preserves parallel lines, and since there always  is  a `nice' orthogonal projection from $d\geq 3$ dimensions into the plane, the best bounds on the number of slopes are obtained in two dimensions. Here a projection is `nice', if no vertex-vertex or vertex-edge occlusions occur; see \citep{EHW97, BGRT95, HW98}. Thus multi-dimensional drawings with few slopes are only interesting if the vertices are restricted to not all lie in a single plane. Under this assumption, \citet{PPS-SoCG04} proved that the minimum number of slopes determined by $n$ points in $\mathbb{R}^3$ is (exactly) $2n-5$ if $n$ is odd, and at least $2n-7$ if $n$ is even. Earlier, \citet{PPS-JCTA04} proved that under the additional assumption that no three points are collinear (which is needed for a drawing of $K_n$), the minimum number of slopes is (exactly) $2n-2$ if $n$ is odd and $2n-3$ if $n$ is even. These proofs are based on generalisations of the above-mentioned result of \citet{Ungar-JCTA82}. In related work, \citet{OP-JCTA04} studied the minimum number of edge slopes in a $d$-dimensional convex polytope. 






\mySection{Graphs of Bounded Degree}{Degree}

Here we prove the following theorem, which was introduced in the introduction.

\begin{theorem}
\thmlabel{DegreeLowerBound}
For all $\Delta\geq5$ and $\varepsilon>0$, for all sufficiently large $n>n(\Delta,\varepsilon)$, there exists a $\Delta$-regular $n$-vertex graph $G$ with slope-number 
\begin{equation*}
\sn{G}\;>\;n^{1-\frac{8+\varepsilon}{\Delta+4}}.
\end{equation*}
\end{theorem}

\begin{proof}
In this proof, $c$ is an positive (absolute) constant that might change from one line to the next. We proceed as in the proof by \citet{BMW-EJC06}. The idea is to show that there are more $\Delta$-regular graphs than $\Delta$-regular graphs with slope-number $k$, for an appropriately chosen $k$. For ease of counting we work with labelled graphs. 

Let \GGG\ be the set of labelled $\Delta$-regular $n$-vertex graphs. The first asymptotic bounds on $|\GGG|$ were independently obtained by \citet{BC-JCTA78} and \citet{Wormald78}. Based on a further refinement by \citet{McKay-AC85}, \citet{BMW-EJC06} proved that 
\begin{equation}
\eqnlabel{NumberRegular}
|\GGG|\geq\BFRAC{n}{3\Delta}^{\Delta n/2}\text{ for all }n\geq c\Delta. 
\end{equation}
The key contribution of \citet{BMW-EJC06} was to show that the number of labelled $n$-vertex $m$-edge graphs with slope-number at most $k$ is at most
\begin{equation}
\eqnlabel{NumberSN}
\BFRAC{50n^2(k+1)}{2n+k}^{2n+k}\binom{k(n-1)}{m}.
\end{equation}

Suppose, on the contrary, that for some $\Delta\geq5$, for some $\varepsilon>0$, and for some $n$, every $\Delta$-regular $n$-vertex graph has slope-number at most 
\begin{equation*}
k:=n^{1-\frac{8+\varepsilon}{\Delta+4}}.
\end{equation*}
We now derive a contradiction for all sufficiently large $n>n(\Delta,\varepsilon)$. By \eqnref{NumberRegular} and \eqnref{NumberSN}, \begin{equation*}
\BFRAC{n}{3\Delta}^{\Delta n/2}
\;\leq\;|\GGG|
\;\leq\;
\BFRAC{50n^2(k+1)}{2n+k}^{2n+k}\binom{k(n-1)}{\Delta n/2}
\;<\;
(ckn)^{2n+k}\binom{kn}{\Delta n/2}.
\end{equation*}
Since $\binom{a}{b}\leq\bfrac{\e\,a}{b}^b$, 
\begin{equation*}
\BFRAC{n}{3\Delta}^{\Delta n/2}
\;<\;
(ckn)^{2n+k}\BFRAC{ck}{\Delta}^{\Delta n/2}.
\end{equation*}
Hence 
\begin{equation*}
n^{4\Delta n}
\;<\;
(ckn)^{16n+8k}(ck)^{4\Delta n}.
\end{equation*}
Observe that $8k<\varepsilon n$ for all large $n>n(\Delta,\varepsilon)$. Thus
\begin{equation*}
n^{4\Delta}<(ckn)^{16+\varepsilon}(ck)^{4\Delta}.
\end{equation*}
That is,
\begin{equation*}
n^{4\Delta-16-\varepsilon}<c^{4\Delta+16+\varepsilon}\,k^{4\Delta+16+\varepsilon}.
\end{equation*}
Since $c^{4\Delta+16+\varepsilon}<n^{2\varepsilon}$ for all large $n>n(\Delta,\varepsilon)$, 
\begin{equation*}
n^{4\Delta-16-3\varepsilon}<k^{4\Delta+16+\varepsilon}.
\end{equation*}
That is,
\begin{equation*}
k
\;>\;n^{\frac{4\Delta-16-3\varepsilon}{4\Delta+16+\varepsilon}}
\;=\;n^{1-\frac{32+4\varepsilon}{4\Delta+16+\varepsilon}}
\;>\;n^{1-\frac{8+\varepsilon}{\Delta+4}},
\end{equation*}
which is the desired contradiction. Therefore for all sufficiently large $n>n(\Delta,\varepsilon)$, there exists a $\Delta$-regular $n$-vertex graph $G$ with $\sn{G}>k$.
\end{proof}

The following open problem remains unsolved.

\begin{open}
Does every graph with maximum degree at most $4$ have bounded slope-number? Note that \citet{DEK-SoCG04} proved that such graphs have geometric thickness at most $2$.
\end{open}

Another interesting problem is to determine the best possible bounds on the slope-number of graphs with bounded degree.

\begin{open}
Does every $n$-vertex graph with bounded degree have $o(n)$ slope-number?
\end{open}

\mySection{Complete Multipartite Graphs}{Multipartite}

We start this section by considering the slope-number of the complete graph $K_n$ on $n$ vertices. Consider a drawing of a graph $G$ on a regular $n$-gon with vertex ordering $(v_1,v_2,\dots,v_n)$.  \citet{Scott70} observed that the number of slopes  is 
\begin{equation}
\eqnlabel{RegularGon}
|\{(i+j)\bmod{n}:v_iv_j\in E(G)\}|\enspace.
\end{equation}
Thus for $K_n$, drawn on a regular $n$-gon, the number of slopes is $n$, as illustrated in \figref{CompleteGraph}. Thus $\sn{K_n}\leq\csn{K_n}\leq n$. To see that this construction is optimal, let $u,v,w$ be three consecutive vertices on the convex hull of an arbitrary drawing of $K_n$. \citet{Jamison-DM86} observed that the $n-1$ edges incident to $v$ and the edge $uw$ have distinct slopes\footnote{More generally, \citet{Jamison-DM86} proved that if a drawing of $K_n$ has $k$ vertices on the convex hull then the number of slopes is at least $k(n-2)/(k-2)$, and that every drawing of $K_n$ with exactly $n$ slopes is affinely equivalent to a regular $n$-gon. Note that \citet{WadeChu-CJ94} independently proved that $\sn{K_n}=n$, and also presented an algorithm to test if $K_n$ can be drawn using a given set of slopes.}. Thus:

\begin{proposition}[\citep{Jamison-DM86}]
$\csn{K_n}=\sn{K_n}=n$.
\end{proposition}

\Figure{CompleteGraph}{\includegraphics{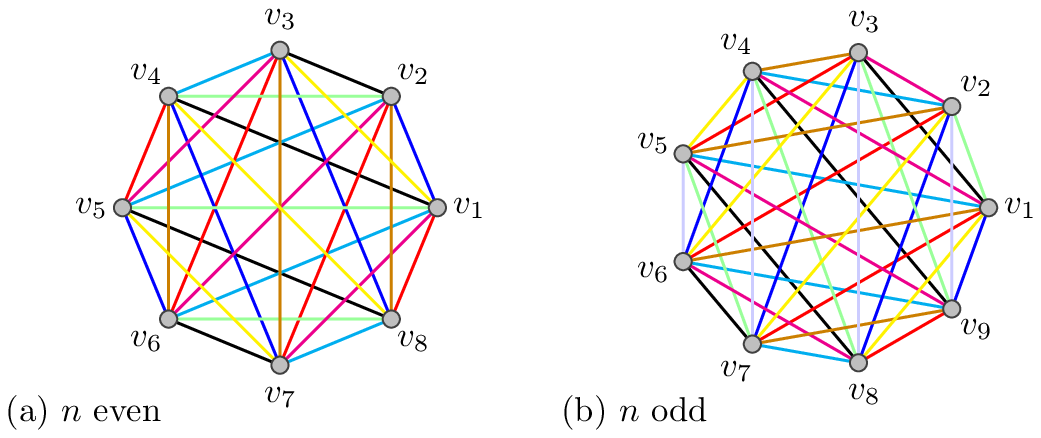}}{Drawings of $K_n$ with $n$ slopes.}

For $k\geq2$, the \emph{complete $k$-partite graph} $K_{n_1,n_2,\dots,n_k}$ has vertex set $V(G):=\{v_{i,j}:1\leq i\leq k,1\leq j\leq n_i\}$ and edge set 
$E(G)=\{v_{i,p}v_{j,q}:1\leq i<j\leq k,1\leq p\leq n_i,1\leq q\leq n_j\}$. 

The slope-number of the balanced complete bipartite graph is easily determined.

\begin{proposition} 
\proplabel{BalancedCompleteBipartiteGraph}
$\sn{K_{n,n}}=\csn{K_{n,n}}=n$.
\end{proposition}

\begin{proof}  
Since $K_{n,n}$ is $n$-regular, $\sn{K_{n,n}}\geq n$ by Equation~(\ref{eqn:LowerBoundSlopes}b). For the upper bound, position the vertices of $K_{n,n}$ on a regular $2n$-gon $(v_1,v_2,\dots,v_{2n})$, alternating between the colour classes, as illustrated in \figref{CompleteBipartiteGraph}. Thus $v_iv_j$ is an edge if and only if $i+j$ is odd. By \eqnref{RegularGon}, the number of slopes is  $|\{(i+j)\bmod{2n}:1\leq i<j\leq 2n,i+j\text{ is odd}\}|=n$.   
\end{proof}

\Figure{CompleteBipartiteGraph}{\includegraphics{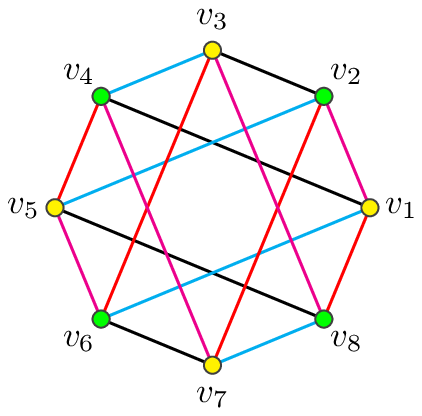}}{Drawing of $K_{4,4}$ with $4$ slopes.}

\propref{BalancedCompleteBipartiteGraph} implies that
$\sn{K_{a,b}}\leq\csn{K_{a,b}}\leq\max\{a,b\}$. In fact,  by Equation~(\ref{eqn:LowerBoundSlopes}c), $\csn{K_{a,b}}\geq\Delta(K_{a,b})=\max\{a,b\}$. Thus 
$\csn{K_{a,b}}=\max\{a,b\}$. Determining $\sn{K_{a,b}}$ is more challenging. We have the following bounds.

\begin{theorem}
\thmlabel{ConvexCompleteBipartiteGraph} 
For all $a\leq b$, 
$\half(a+b-1)\leq \sn{K_{a,b}}\leq\min\{b,\ceil{\frac{b}{2}}+a-1\}$.
\end{theorem}

\begin{proof} 
That $\sn{K_{a,b}}\leq b$ follows from \propref{BalancedCompleteBipartiteGraph}.
 
Now we prove the upper bound, $\sn{K_{a,b}}\leq\ceil{\frac{b}{2}}+a-1$. Without loss of generality $b$ is even. Suppose $V(K_{a,b})=\{v_1,v_2,\dots,v_a\}\cup\{u_1,u_2,\dots,u_{\frac{b}{2}}\}\cup\{w_1,w_2,\dots,w_{\frac{b}{2}}\}$, and $E(K_{a,b})=\{v_iu_j,v_iw_j:1\leq i\leq a,1\leq j\leq \frac{b}{2}\}$.  Position each vertex $u_j$ at $(j,1)$; position each vertex $v_i$ at $(\frac{b}{2}+i,0)$; and position each vertex $w_j$ at $(\frac{b}{2}+a+j,-1)$. Then every edge is parallel with one of the $\frac{b}{2}+a-1$  edges $\{v_1u_j:1\leq j\leq \frac{b}{2}\}\cup\{u_1v_i:2\leq i\leq a\}$, as illustrated in \figref{NewCompleteBipartiteGraph}.

\Figure{NewCompleteBipartiteGraph}{\includegraphics{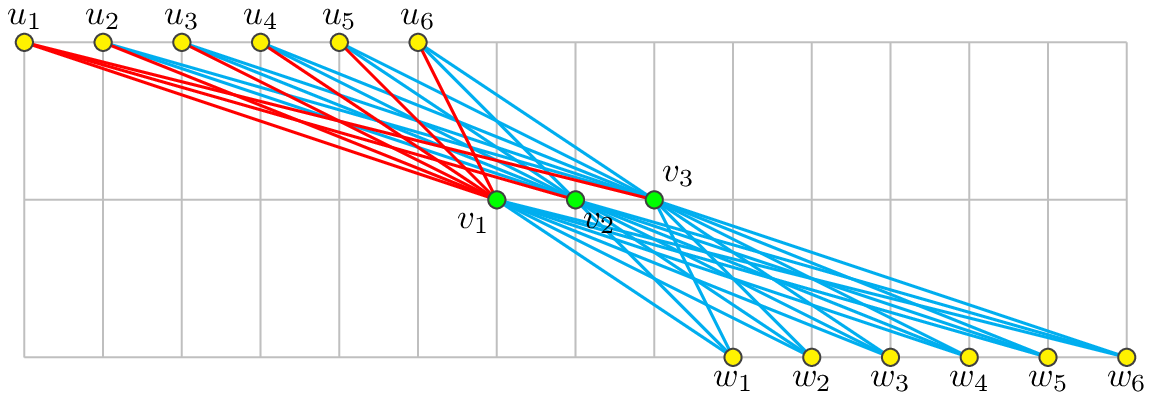}}{
Drawing of $K_{3,12}$ with $8$ slopes (highlighted).}

Now we prove the lower bound (which is due to an anonymous referee). Let $A$ and $B$ be the two colour classes of $K_{a,b}$ where $|A|=a$ and $|B|=b$. Given a drawing of $K_{a,b}$, rotate it so that no two vertices are horizontal. Let $L$ be a horizontal line that intersects no vertex, and has at least $\floor{\half(a+b)}$ vertices above and below $L$. Let $a_1$ and $a_2$ be the number of vertices in $A$ respectively above and below $L$. Let $b_1$ and $b_2$ be the number of vertices in $B$ respectively above and below $L$. 
Thus $a_1+b_1\geq\floor{\half(a+b)}$ and $a_2+b_2\geq\floor{\half(a+b)}$.
Since $(a_1+b_2)+(a_2+b_1)=a+b$, without loss of generality, $a_1+b_2\geq\ceil{\half(a+b)}$. 

We claim that $a_1>0$ and $b_2>0$. Suppose on the contrary that $a_1=0$. Thus $b=b_1+b_2\geq\floor{\half(a+b)}+\ceil{\half(a+b)}=a+b$, implying $a=0$, which is a contradiction. Thus $a_1>0$, and similarly, $b_2>0$. Consider the drawing of $K_{a_1,b_2}$ induced by the $a_1$ vertices in $A$ above $L$ and the $b_2$ vertices in $B$ below $L$. Every edge of $K_{a_1,b_2}$ crosses $L$ and there is some edge in $K_{a_1,b_2}$. Let $vw$ be the leftmost edge of $K_{a_1,b_2}$ crossing $L$. Then the $a_1+b_2-1$ edges of $K_{a_1,b_2}$ incident to $v$ or $w$ all have distinct slopes, as illustrated in \figref{BipartiteLowerBound}.
\end{proof}


\Figure{BipartiteLowerBound}{\includegraphics{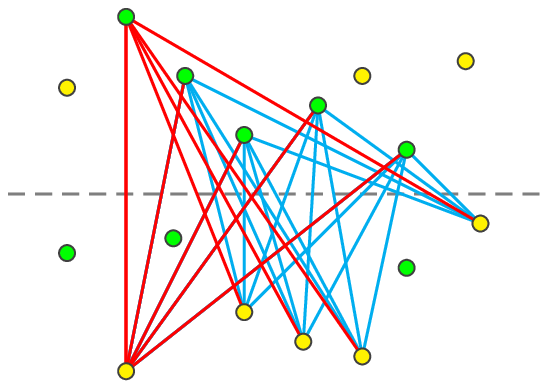}}{Finding a large separated subgraph in $K_{a,b}$.}

Closing the gap in the bounds in \thmref{ConvexCompleteBipartiteGraph} remains an interesting open problem.

\begin{open}
\openlabel{CompleteBipartite}
What is the slope-number \sn{K_{a,b}} of the complete bipartite graph $K_{a,b}$?
\end{open}

Now consider the general case of a complete $k$-partite graph $G$. 
Say $G$ has $n$ vertices. Since $\csn{G}\leq n$ and $\Delta(G)\geq\frac{k-1}{k}n$, we have $\sn{G}\leq\csn{G}\leq\frac{k}{k-1}\Delta(G)$.

\begin{open}
\openlabel{CompleteMultipartiteGraph}
Does every complete multipartite graph $G$ with maximum degree $\Delta$ have a (convex) drawing with at most $\Delta+\text{o}(\Delta)$ slopes?
\end{open}

We have the following partial solution to \openref{CompleteMultipartiteGraph}.

\begin{proposition}
\proplabel{PowerTwoComplete}
Given integers $p\geq0$ and $k\geq2$, where $k-1$ is a power of two, let $G$ be the complete $k$-partite graph $K_{2^p,2^p,2^{p+1},\ldots,2^{p+1}}$. Then $\sn{G}\leq\csn{G}=\Delta(G)$.
\end{proposition}

\begin{proof}
Equation~(\ref{eqn:LowerBoundSlopes}c) implies that $\csn{G}\geq \Delta(G)$. We now prove the upper bound. 

Let $n:=(k-1)2^{p+1}$ be the number of vertices in $G$. Note that $n$ is a power of two, and $\Delta(G)=n-2^p$. In what follows $a\equiv b$ means that $a\equiv b\pmod{n/2^p}$, and $a\equiv\pm b$ means that $a\equiv b$ or $a\equiv -b$. 
For all $0\leq i\leq k-1$, let $P_i=\{j\in V(G)\,:\,i\equiv\pm j\}$. 
Let $V(G):=\{0,1,\dots,n-1\}$. Below we prove that $\{P_0,P_1,\ldots,P_{k-1}\}$ is a partition of $V(G)$ with $|P_0|=|P_{k-1}|=2^p$, and $|P_i|=2^{p+1}$ for all $1\leq i\leq k-2$.  Thus $\{P_0,P_1,\ldots,P_{k-1}\}$ defines a valid assignment of the vertices to the colour classes. To obtain the drawing of $G$, place the vertices in numerical order on the vertices of a regular $n$-gon. 

For each vertex $j\in V(G)$, let $j':=j\bmod{n/2^p}$. If $0\leq j'\leq n/2^{p+1}$, then $j\in P_{j'}$. Otherwise, $n/2^{p+1} < j' < n/2^p$, and $j\in P_{n/2^p-j'}$. Thus, each vertex belongs to at least one $P_i$. Suppose that $j\in P_i\cap P_h$. Thus $i\equiv\pm j$ and $h\equiv\pm j$, implying $i\equiv\pm h$. Since $0\leq i\leq n/2^{p+1}$, we have $h=i$. Thus, each vertex belongs to exactly one $P_i$, and  $\{P_0,P_1,\ldots,P_{k-1}\}$ is a partition of  $V(G)$. The set $P_0$ has size $2^p$ because   it is the set of all multiples of $n/2^p$ in $\{0,1,\dots,n-1\}$. Similarly, $P_{k-1}$ has size $2^p$ because  it is the set of all odd multiples of $n/2^{p+1}$ in $\{0,1,\dots,n-1\}$. The remainder of the $P_i$'s have the same size, $2^{p+1}$, by symmetry. 

To prove that the number of slopes $|\{(i+j)\bmod{n}:ij\in E(G)\}|=n-2^p$, by \eqnref{RegularGon}, it suffices to prove that $\mequiv{i+j}{0}{n/2^p}$ implies $ij\not\in E(G)$. Suppose that $i\in P_h$. Thus $h+i\equiv 0$ or $\mequiv{h-i}{0}{n/2^p}$. In the first case, we have $\mequiv{h+i}{i+j}{n/2^p}$, implying $\mequiv{h-j}{0}{n/2^p}$. In the second case, we have $\mequiv{h-i+(i+j)}{0}{n/2^p}$, implying  $\mequiv{h+j}{0}{n/2^p}$. In both cases $j\in P_h$, implying $ij\not\in E(G)$. 
\end{proof}

\begin{corollary}
\corlabel{PowerTwoComplete}
Given integers $p\geq0$, $q\leq 2^p$, and $k\geq2$, where $k-1$ is a power of two, let $G$ be the complete $k$-partite graph $K_{q,2^p,2^{p+1},\ldots,2^{p+1}}$. Then $\csn{G}=\Delta(G)$.
\end{corollary}

\begin{proof}
Let $G'$ be the complete $k$-partite graph $K_{2^p,2^p,2^{p+1},\ldots,2^{p+1}}$. Then $G$ is a subgraph of  $G'$, and $\Delta(G)=\Delta(G')= (k-2)2^{p+1}+2^p$.  The result follows from \propref{PowerTwoComplete}. \end{proof}


\mySection{General Graphs}{GeneralSlope}

While \thmref{DegreeLowerBound} proves that there exist graphs of bounded degree with unbounded slope-number, in this section, we prove that the slope-number of various classes of graphs is bounded by a function of the maximum degree. For graphs of bounded degree and bounded treewidth we prove a \Oh{\log n} bound on the slope-number.

Our results are based on the following structure. Let $H$ be a (\emph{host}) graph. The vertices of $H$ are called \emph{nodes}. An \emph{$H$-partition} of a graph $G$ is a function $f:V(G)\rightarrow V(H)$ such that for every edge $vw\in E(G)$ we have $f(v)=f(w)$ or $f(v)f(w)\in E(H)$. In the latter case, we say $vw$ is \emph{mapped} to the edge $f(v)f(w)$. The \emph{width} of $f$ is the maximum of $|f^{-1}(x)|$, taken over all nodes $x\in V(H)$, where $f^{-1}(x):=\{v\in V(G):f(v)=x\}$. The following general result describes how to produce a drawing of a graph $G$ given an $H$-partition of $G$ and a drawing of $H$.

\begin{theorem}
\thmlabel{BlowUp}
Let $D$ be a drawing of a graph $H$ with $s$ distinct slopes and $\ell$ distinct edge lengths. Let $t:=|\{(\textup{slope}(e),\textup{length}(e)):e\in E(D)\}|$ \paran{which is at most $s\ell$}. Let $G$ be a graph with an $H$-partition of width $k$. Then $G$ has a drawing with at most $k+s+t(k^2-k)\leq k+s+s\ell(k^2-k)$ distinct slopes \paran{and at most $\floor{\frac{k}{2}}+\ell+t(k^2-k)\leq \floor{\frac{k}{2}}+\ell+s\ell(k^2-k)$ distinct edge lengths}. 
\end{theorem}

\begin{proof} 
The general approach is to scale $D$ appropriately, and then replace each node of $H$ by a copy of the drawing of $K_k$ on a regular $k$-gon (described in \secref{Multipartite}). The only difficulty is to scale $D$ so that we obtain a valid drawing of $G$.

Observe that $|\phi_1-\phi_2|$ is the size of the minimum angle formed by lines of slope $\phi_1$ and $\phi_2$. Let $\{\theta_1,\theta_2,\dots,\theta_s\}$ be the set of slopes of the edges of $D$. Rotate the drawing of $K_k$ on a regular $k$-gon so that, if $\{\beta_1,\beta_2,\dots,\beta_k\}$ is the set of slopes of the edges of $K_k$, then $|\theta_i-\beta_j|>0$ for all $i$ and $j$. Let $\varepsilon:=\min_{i,j}\{|\theta_i-\beta_j|\}$.

Replace each node $x$ in $D$ by a disc $B_x$ of uniform radius $r$ centred at
$x$, where $r$ is chosen small enough so that: (1)  $B_x\cap B_y=\emptyset$ for
all distinct nodes $x$ and $y$ in $D$; and (2) for every edge $xy\in E(H)$ with
slope $\theta_i$, every segment with endpoints in $B_x$ and $B_y$ and with
slope $\phi$ intersects no other $B_z$,  and $|\phi-\theta_i|<\varepsilon$. 
Position a regular $k$-gon on each $B_x$ (using the orientation determined
above), and position the vertices $f^{-1}(x)$ of $G$ at its vertices. Since
$|\theta_i-\beta_j|\geq\varepsilon$, the slope of any edge $vw$ of $G$ that
is mapped to $xy$ does not equal any $\beta_j$. Hence $vw$ does not pass
through any other vertex of $G$. 

Each copy of $K_k$ contributes the same $k$ slopes to the drawing of $G$.
For each edge $xy\in E(H)$, for all $1\leq i\leq k$, the edge of $G$ from the
$i$-th vertex on $B_x$ to the $i$-th vertex on $B_y$ (if it exists) has the same slope as the edge $xy$ in $D$. Thus these edges contribute $s$ slopes to the drawing of $G$. Consider two edges $e_1$ and $e_2$ of $H$ that have the same slope and the same length in $D$ (of the $t$ possibilities). The edges of $G$ that are mapped to $e_1$ use the same set of slopes as the edges of $G$ that are mapped to $e_2$. There are at most $k^2-k$ edges of $G$ that are mapped to a single edge of $H$ and were not counted above. Thus in total we have at most $k+s+t(k^2-k)$ slopes, as illustrated in \figref{BlowUp}.

Each copy of $K_k$ contributes the same $\ceil{\frac{k}{2}}$ distinct edge
lengths. This, along with analogous arguments to those presented above, gives
an upper bound of $\floor{\frac{k}{2}}+\ell+t(k^2-k)$ on the number of distinct edge lengths.
\end{proof}

\Figure{BlowUp}{\includegraphics{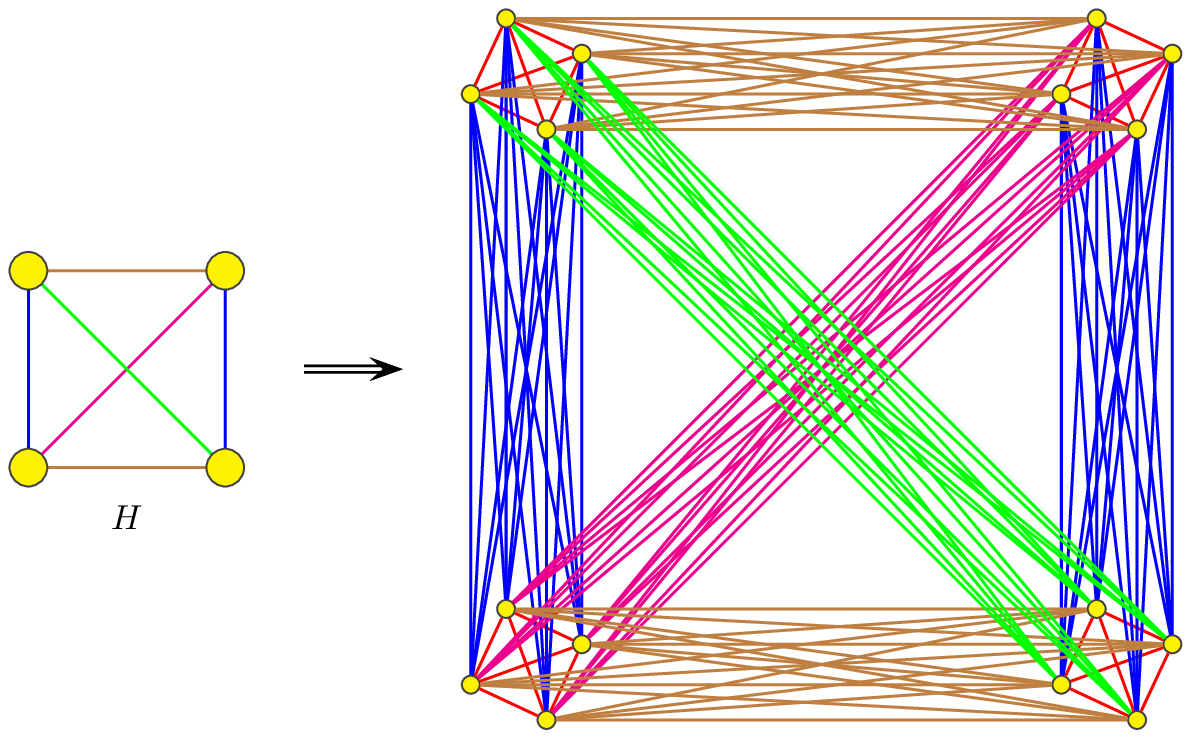}}{Illustration of the construction in \thmref{BlowUp} with $H=K_4$, $s=4$, $\ell=2$, and $k=4$.}

\mySubSection{Drawings Based on Paths}{Paths}

\thmref{BlowUp} suggests using host graphs that have drawings with few
slopes and few edge lengths. Thus a path is a natural choice for a host graph, since it has a drawing with one slope and one edge length. The \emph{path-partition-width} of a graph $G$, denoted by $\ppw{G}$, is the minimum integer $k$ such that $G$ has a $P$-partition of width $k$, for some path $P$. \thmref{BlowUp} with $r=s=\ell=1$ implies:

\begin{corollary}
\corlabel{PPWandSlopes}
Every graph $G$ has a drawing with $\ppw{G}^2+1$ slopes.\qed 
\end{corollary}

As indicated by the following lemma, path-partition-width is closely related to the classical graph parameter bandwidth\footnote{\citet{Bodlaender-IC90} found essentially the same relation in the context of emulations of networks.}. The \emph{width} of a vertex ordering $(v_1,v_2,\dots,v_n)$ of a graph $G$ is the maximum of $|i-j|$, taken over  all edges $v_iv_j\in E(G)$. The \emph{bandwidth} of $G$, denoted by \bw{G}, is the minimum width of a vertex ordering of $G$.

\begin{lemma}
\lemlabel{PPWandBW}
For every graph $G$, $\half(\bw{G}+1)\leq\ppw{G}\leq\bw{G}$.
\end{lemma}

\begin{proof} Let $(v_1,v_2,\dots,v_n)$ be a vertex ordering of $G$ with width
$b=\bw{G}$. For all $0\leq i\leq\floor{n/b}$,  let
$B_i:=\{v_{ib+1},v_{ib+2},\dots,v_{ib+b}\}$. Then $(B_0,B_1,\dots,B_{\floor{n/b}})$ defines a path-partition  of $G$ with width
$b$.  Thus $\ppw{G}\leq\bw{G}$.

Now suppose $(B_1,B_2,\dots,B_m)$ is a path-partition of $G$ with width
$k=\ppw{G}$. Let $(v_1,v_2,\dots,v_n)$ be a vertex ordering of $G$ such that
$i<j$ whenever $v_i\in B_p$ and $v_j\in B_q$ and $p<q$. For every edge
$v_iv_j\in E(G)$ with $v_i\in B_p$ and $v_j\in B_q$, we have $|p-q|\leq1$. Thus
$|i-j|\leq 2k-1$. Hence the width of $(v_1,v_2,\dots,v_n)$ is at most $2k-1$.
Therefore $\bw{G}\leq 2\ppw{G}-1$.
\end{proof}

\corref{PPWandSlopes} and \lemref{PPWandBW} imply that every graph $G$ has a
drawing with $\bw{G}^2+1$ slopes. This bound can be tweaked as follows.

\begin{theorem}
\thmlabel{Bandwidth}
Every graph $G$ has slope-number $\sn{G}\leq\half\bw{G}\,(\bw{G}+1)+1$.
\end{theorem}

\begin{proof}  
Let $G[B_i,B_{i+1}]$ be the bipartite subgraph of $G$ with vertex set $B_i\cup B_{i+1}$ and edge set $\{vw\in E(G):v\in B_i,w\in B_{i+1}\}$ from the proof of \lemref{PPWandBW}. Observe that in the construction of the path-partition in \lemref{PPWandBW}, the edges of each $G[B_i,B_{i+1}]$ are a subset of
$\{\{v_{ib+j},v_{(i+1)b+\ell}\}:1\leq j\leq b,1\leq \ell\leq j \}$.  If we
consistently assign the vertices in each $B_i$ to the regular $b$-gon in
\thmref{BlowUp}, then each $G[B_i,B_{i+1}]$ will use the same set of slopes,
since each $G[B_i, B_{i+1}]$ is a subgraph of the same graph.  The number of
slopes in $G[B_i, B_{i+1}]$ is $1+\sum_{j=1}^b(j-1)$, since  each vertex
$v_j\in B_i$ is incident to $j$ edges with endpoints in $B_{i+1}$, one of which
is horizontal. Thus the total number of slopes in the resulting drawing of $G$
is $b+1+\half(b-1)b=\half b(b+1)+1$. \end{proof}

\Figure{Bandwidth}{\includegraphics{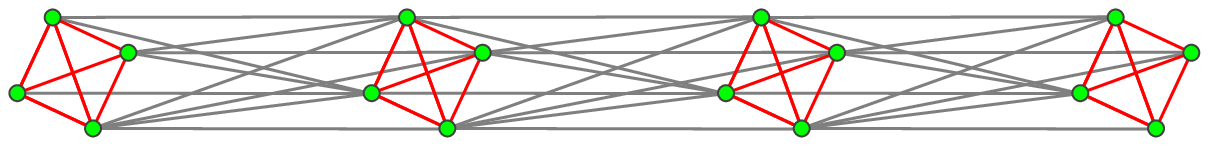}}{Drawing of a graph with
bandwidth $4$ with eleven slopes.}

The following examples of \thmref{Bandwidth} are corollaries of results by
\citet{FominGolovach-DAM03} and \citet{Wood-IntersectionGraphOrderings} that
bound bandwidth in terms of maximum degree.

\begin{itemize}

\item Every interval graph $G$  has $\bw{G}\leq\Delta(G)$
\citep{FominGolovach-DAM03,Wood-IntersectionGraphOrderings}, and thus has a
drawing with at most  $\half\Delta(G)\,(\Delta(G)+1)+1$ slopes.

\item Every cocomparability graph $G$ has $\bw{G}\leq2\Delta(G)-1$ \citep{Wood-IntersectionGraphOrderings}, and thus
has a drawing with at most $\Delta(G)\,(2\Delta(G)-1)+1$ slopes.

\item Every AT-free graph $G$ has $\bw{G}\leq3\Delta(G)$
\citep{Wood-IntersectionGraphOrderings}, and thus has a drawing with at most
$\frac{3}{2}\Delta(G)\,(3\Delta(G)+1)+1$ slopes.

\end{itemize}

\begin{open} 
Does every interval graph $G$ have a drawing with \Oh{\Delta(G)} slopes? 
\end{open}

\mySubSection{Drawings Based on Trees}{Trees}

To obtain bounds on the slope-number of more general graphs, we consider $T$-partitions for some tree $T$. This structure is called a \emph{tree-partition}, and has been extensively studied \citep{Halin91, BodEng-JAlg97, Wood-TreePartitions, Bodlaender-DMTCS99, Halin91, Seese85, DO-JGT95, DO-DM96, Edenbrandt86, DMW-SJC05}. \thmref{BlowUp} motivates the study of drawings of trees with few slopes and few distinct edge lengths.

\begin{theorem}
\thmlabel{TreeSlopesLengths}
Every tree $T$ with pathwidth $k\geq1$ has a plane drawing  with
$\max\{\Delta(T)-1,1\}$ slopes and $2k-1$ distinct edge lengths.
\end{theorem}


The proof of \thmref{TreeSlopesLengths} is loosely based  on an algorithm of
\citet{Suderman-IJCGA04} for drawing trees on layers. We will need the following
lemma\footnote{In fact, \lemref{MainPath} can be viewed as the basis for an alternative definition of the pathwidth of a forest. In particular, the pathwidth of $K_1$ equals $0$, the pathwidth of a forest $F$ equals the maximum pathwidth of a connected component of $F$, and the pathwidth of a tree $T$ equals the minimum $k$ such that there exists a path $P$ of $T$ and the pathwidth of $T\setminus V(P)$ is at most $k-1$.}. 

\begin{lemma}[\citep{Suderman-IJCGA04}]
\lemlabel{MainPath} 
Every tree $T$ has a path $P$ such that $T\setminus V(P)$ has smaller pathwidth
than $T$, and the endpoints of $P$ are leaves of $T$. 
\end{lemma}


A path $P$ satisfying \lemref{MainPath} is called a \emph{backbone} of $T$.

\begin{proof}[Proof of \thmref{TreeSlopesLengths}]  
We refer to $T$ as $T_0$. Let $n_0$ be the number of vertices in $T_0$, and let
$\Delta_0=\Delta(T_0)$. The result holds trivially for $\Delta_0\leq 2$. Now assume that $\Delta_0\geq3$. Let  $S$ be  the set of slopes
\begin{equation*}
S\;:=\;
\SET{\frac{\pi}{2}\BRACKET{1+\frac{i}{\Delta_0-2}}\,:\,0\leq i\leq \Delta_0-2}
\enspace.
\end{equation*}

\noindent We proceed by induction on $n$ with the hypothesis:  ``There is a
real number $\ell=\ell(n_0,\Delta_0)$, such that  for every tree $T$ with
$n\leq n_0$ vertices, maximum degree at most $\Delta_0$, and pathwidth $k\geq1$, 
and for every vertex $r$ of $T$ with degree less than $\Delta_0$,  $T$ has a
plane drawing $D$ in which:
\begin{itemize}
\item $r$ is at the top of $D$ (that is, no point in $D$ has greater
Y-coordinate than $r$),
\item every edge of $T$ has slope in $S$,
\item every edge of $T$ has length in $\{1,\ell,\dots,\ell^{2k-1}\}$, and
\item if $r$ is contained in some backbone of $T$, then 
every edge of $T$ has length in $\{1,\ell,\dots,\ell^{2k-2}\}$.''
\end{itemize}

The result follows from the induction hypothesis, since we can take $r$ to be
the endpoint of a backbone of $T_0$, in which case $\deg(r)=1<\Delta_0$, and
thus every edge of $T_0$ has length in $\{1,\ell,\dots,\ell^{2k-2}\}$.

The base case with $n=1$ is trivial. Now suppose that the hypothesis is true
for trees on less than $n$ vertices, and we are given a tree $T$ with $n$
vertices and pathwidth $k$, and $r$ is a vertex of $T$ with degree less than
$\Delta_0$.

If $r$ is contained in some backbone $B$ of $T$, then let $P:=B$. Otherwise,  let $P$ be a path from $r$ to an endpoint of a backbone $B$ of $T$. Note that
$P$ has at least one edge.  As illustrated in \figref{TreeLengthsSlopes}, draw
$P$ horizontally with  unit-length edges. Every vertex in $P$ has at most
$\Delta_0-2$ neighbours in $T\setminus V(P)$, since $r$ has degree less than
$\Delta_0$ and the endpoints of a backbone are leaves. At each vertex $x\in P$,
the children $\{y_0,y_1,\dots,y_{\Delta_0-3}\}$ of $x$ are positioned below $P$
and on the unit-circle centred at $x$, so that each edge $xy_j$ has slope
$\frac{\pi}{2}(1+j/(\Delta_0-2))\in S$. 

\Figure{TreeLengthsSlopes}{\includegraphics{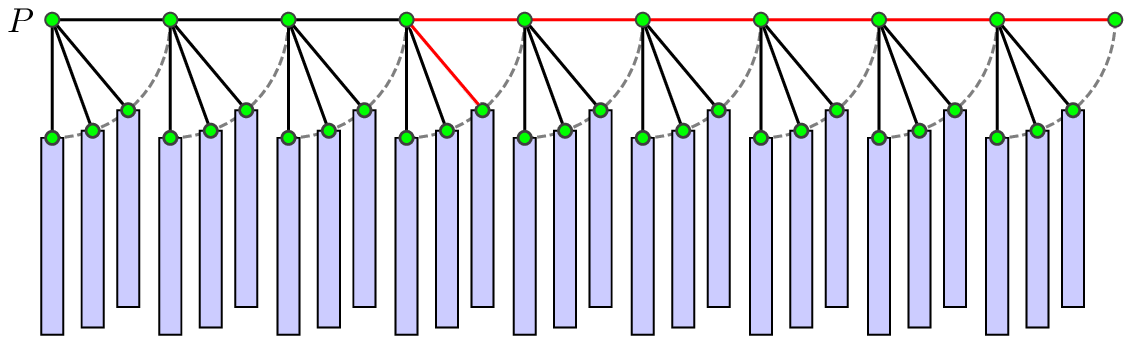}}{Drawing of $T$
with few slopes and few edge lengths.}

Every connected component $T'$ of $T\setminus V(P)$ is a tree  rooted at some
vertex $r'$ adjacent to a vertex in $P$. By the above layout procedure, $r'$
has already been positioned in the drawing of $T$. If $T'$ is a single vertex,
then we no longer need to consider this $T'$. 

We consider two types of subtrees $T'$, depending on whether the pathwidth of
$T'$ is less than $k$. Suppose that the pathwidth of $T'$ is $k$ (it cannot be
more). Then $T'\cap B\ne\emptyset$ since $B$ is a backbone of $T$. Thus $T'\cap
B$ is a backbone of $T'$ containing $r'$. Thus we can apply the stronger
induction hypothesis in this case.

Every $T'$ has fewer vertices than $T$, and every $r'$ has degree less than
$\Delta_0$ in $T'$. Thus by induction,  every    $T'$ has a drawing with $r'$ at
the top, and every edge of $T'$ has slope in $S$.  Furthermore, if the
pathwidth of $T'$ is less than $k$, then  every edge of $T'$ has length in
$\{1,\ell,\dots,\ell^{2k-3}\}$. Otherwise $r'$ is in a backbone of $T'$,
and every edge of $T'$ has length in $\{1,\ell,\dots,\ell^{2k-2}\}$.

There exists a scale factor $\ell<1$, depending only on $n_0$ and $\Delta_0$,  so that by scaling the drawings of every $T'$ by $\ell$,  the widths of the
drawings are small enough so that there is no crossings when the drawings are 
positioned with each $r'$ at its already chosen location. (Note that $\ell$ is
the same value at every level of the induction.)\ Scaling preserves the slopes
of the edges. An edge in any $T'$ that had length $\ell^i$ before scaling, now
has length $\ell^{i+1}$. 

Case 1. $r$ is contained in some backbone $B$ of $T$: By construction, $P=B$.
So every $T'$ has pathwidth at most $k-1$, and thus every edge  of $T'$ has
length in $\{\ell^1,\ell^2,\dots,\ell^{2k-2}\}$. All the other edges of $T$
have unit-length. Thus we have a plane drawing of $T$ with edge lengths 
$\{1,\ell,\dots,\ell^{2k-2}\}$, as claimed.

Case 2. $r$ is not contained in any backbone of $T$:   Every edge in every $T'$
has length in $\{\ell^1,\ell^2,\dots,\ell^{2k-1}\}$. All the other edges of $T$
have unit-length. Thus we have a plane drawing of $T$ with edge lengths
$\{1,\ell,\dots,\ell^{2k-1}\}$, as claimed. \end{proof}


\begin{theorem}
Let $G$ be a graph with $n$ vertices, maximum degree $\Delta\geq1$, and treewidth
$k\geq1$. Then $G$ has a drawing with \Oh{k^3\Delta^4\log n} slopes.
\end{theorem}

\begin{proof} 
\citet{Wood-TreePartitions} proved that $G$ has a $T$-partition of width at most $w:=2\big(k+1\big)\big(9\,\Delta-1\big)$ for some forest $T$. (The proof is a minor improvement to a similar result by an anonymous referee of the paper by \citet{DO-JGT95}.)\ For each node $x\in V(T)$, there are at most $w\Delta$
edges of $G$ incident to vertices mapped to $x$. Hence we can assume that $T$
is a forest with maximum degree at most $w\Delta$, as otherwise there is an
edge of $T$ with no edge of $G$ mapped to it, in which case the edge of $T$ can
be deleted. Similarly, $T$ has at most $n$ vertices. \citet{Scheffler89} proved
that $T$ has pathwidth at most $\log(2n+1)$; see \citep{Bodlaender-TCS98}. By
\thmref{TreeSlopesLengths}, $T$ has a drawing with at most $w\Delta-1$ slopes
and at most $2\log(2n+1)-1$ distinct edge lengths. By \thmref{BlowUp},  $G$ has
a drawing in which the number of slopes is at most  
$w(w\Delta-1)(2\log(2n+1)-1)(w-1)+(w\Delta-1)+w \in \Oh{w^3\Delta\log 
n}\subseteq\Oh{k^3\Delta^4\log n}$. 
\end{proof}

\begin{corollary}
Every $n$-vertex graph with bounded degree and bounded treewidth has a drawing
with \Oh{\log n} slopes.\qed
\end{corollary}

\mySection{\boldmath $1$-Bend Drawings}{Bends}

While \thmref{DegreeLowerBound} proves that some graph with bounded degree has unbounded slope-number, we now show that there is no such graph if we allow bends in the edges. For a graph $G$, let $G'$ be the graph obtained from $G$ by subdividing each edge of $G$; that is, for each edge $e=vw$ of $G$, introduce a new \emph{subdivision} vertex $x_e$ in $G'$, and replace $e$ by the path $vx_ew$. A \emph{1-bend drawing} of $G$ is a drawing $G'$.

\begin{theorem} 
Every graph $G$ has a $1$-bend drawing with $\Delta(G)+1$ slopes.
\end{theorem}

\begin{proof} Let $S$ be a set of $\Delta(G)+1$ distinct slopes. Suppose the
vertices of $G$ have been positioned in the plane. For each vertex $v$ of $G$
and each slope $\ell\in S$, consider there to be a \emph{slope line} through
$v$ with slope $\ell$. Position the vertices of $G$ at distinct points in the
plane so that: (1) each slope line intersects exactly one vertex, and (2) no
three slope lines intersect at a single point, unless all three are the slope
lines of a single vertex. This can be achieved by positioning each vertex in
turn, since at each step, there are finitely many forbidden positions. 

Consider each slope line to be initially \emph{unused}.  Each edge is drawn
with one bend, using one slope line at each of its endpoints, in which case, we
say these slope lines become \emph{used}. Now draw each edge $vw$ of $G$ in
turn. At most $\deg(v)-1$ slope lines at $v$ are used,  and at most $\deg(w)-1$
slope lines at $w$ are used.  Since $|S|\geq\deg(v)+1$ and $|S|\geq\deg(w)+1$,
there are two unused slope lines at $v$, and two unused slope lines at $w$.
Thus there is an unused slope line at $v$ that intersects an unused slope line
at $w$.  Position the bend for $vw$ at their intersection point. 

We now prove that this defines a drawing of $G'$. Suppose on the contrary that
there is an edge $vu$ of $G'$ and a vertex $w$ of $G'$ that intersects $vu$,
and $v\ne w\ne u$. Without loss of generality, $v$ is a vertex of $G$ and $u$
is a subdivision vertex. Since each slope line intersects exactly one vertex of
$G$, $w$ is a subdivision vertex of some edge $w_1w_2$ of $G$. Since edges are
only drawn on unused slope lines, $w_1\ne v$ and $w_2\ne v$. Therefore, the
three slope lines containing the edges $w_1w$, $w_2w$ and $vu$ intersect in one
point, and all three do not belong to the same vertex.  This is a desired
contradiction. \end{proof}


\section*{Acknowledgements}

This research was initiated at the \emph{International Workshop on Fixed Parameter Tractability in Geometry and Games}, organised by Sue Whitesides; Bellairs Research Institute of McGill University, Barbados, February 7--13, 2004. Thanks to all of the participants for creating a stimulating working environment. Special thanks to Mike Fellows for suggesting the topic. Thanks to the anonymous referees for many helpful suggestions, including the proof of the lower bound in \thmref{ConvexCompleteBipartiteGraph}, and for pointing out reference \citep{Bodlaender-IC90}.


\begin{thebibliography}{49}
\providecommand{\natexlab}[1]{#1}
\providecommand{\url}[1]{\texttt{#1}}
\providecommand{\urlprefix}{}
\expandafter\ifx\csname urlstyle\endcsname\relax
  \providecommand{\doi}[1]{doi:\discretionary{}{}{}#1}\else
  \providecommand{\doi}{doi:\discretionary{}{}{}\begingroup
  \urlstyle{rm}\Url}\fi

\bibitem[{Ambrus et~al.(2005)Ambrus, Bar{\'a}t, and Hajnal}]{ABH05}
\textsc{Gergely Ambrus, J{\'a}nos Bar{\'a}t, and P{\'e}ter Hajnal}.
\newblock The slope parameter of graphs.
\newblock Tech. Rep. MAT-2005-07, Department of Mathematics, Technical
  University of Denmark, Lyngby, Denmark, 2005.

\bibitem[{Bar{\'a}t et~al.(2006)Bar{\'a}t, Matou{\v{s}}ek, and
  Wood}]{BMW-EJC06}
\textsc{J{\'a}nos Bar{\'a}t, Ji{\v{r}}{\'i} Matou{\v{s}}ek, and David~R. Wood}.
\newblock Bounded-degree graphs have arbitrarily large geometric thickness.
\newblock \emph{Electron. J. Combin.}, 13(1):R3, 2006.

\bibitem[{Bender and Canfield(1978)}]{BC-JCTA78}
\textsc{Edward~A. Bender and E.~Rodney Canfield}.
\newblock The asymptotic number of labeled graphs with given degree sequences.
\newblock \emph{J. Combin. Theory Ser. A}, 24:296--307, 1978.

\bibitem[{Bodlaender(1990)}]{Bodlaender-IC90}
\textsc{Hans~L. Bodlaender}.
\newblock The complexity of finding uniform emulations on paths and ring
  networks.
\newblock \emph{Inform. and Comput.}, 86(1):87--106, 1990.

\bibitem[{Bodlaender(1998)}]{Bodlaender-TCS98}
\textsc{Hans~L. Bodlaender}.
\newblock A partial $k$-arboretum of graphs with bounded treewidth.
\newblock \emph{Theoret. Comput. Sci.}, 209(1-2):1--45, 1998.

\bibitem[{Bodlaender(1999)}]{Bodlaender-DMTCS99}
\textsc{Hans~L. Bodlaender}.
\newblock A note on domino treewidth.
\newblock \emph{Discrete Math. Theor. Comput. Sci.}, 3(4):141--150, 1999.

\bibitem[{Bodlaender and Engelfriet(1997)}]{BodEng-JAlg97}
\textsc{Hans~L. Bodlaender and Joost Engelfriet}.
\newblock Domino treewidth.
\newblock \emph{J. Algorithms}, 24(1):94--123, 1997.

\bibitem[{Bose et~al.(1996)Bose, G\'{o}mez, Ramos, and Toussaint}]{BGRT95}
\textsc{Prosenjit Bose, Francisco G\'{o}mez, Pedro~A. Ramos, and Godfried~T.
  Toussaint}.
\newblock Drawing nice projections of objects in space.
\newblock In \textsc{Franz~J. Brandenburg}, ed., \emph{Proc. International
  Symp. on Graph Drawing (GD '95)}, vol. 1027 of \emph{Lecture Notes in Comput.
  Sci.}, pp. 52--63. Springer, 1996.

\bibitem[{Czyzowicz(1991)}]{Czy-JCTA91}
\textsc{Jurek Czyzowicz}.
\newblock Lattice diagrams with few slopes.
\newblock \emph{J. Combin. Theory Ser. A}, 56(1):96--108, 1991.

\bibitem[{Czyzowicz et~al.(1990{\natexlab{a}})Czyzowicz, Pelc, and
  Rival}]{CPR-DM90}
\textsc{Jurek Czyzowicz, Andrzej Pelc, and Ivan Rival}.
\newblock Drawing orders with few slopes.
\newblock \emph{Discrete Math.}, 82(3):233--250, 1990{\natexlab{a}}.

\bibitem[{Czyzowicz et~al.(1990{\natexlab{b}})Czyzowicz, Pelc, Rival, and
  Urrutia}]{CPRU-Order90}
\textsc{Jurek Czyzowicz, Andrzej Pelc, Ivan Rival, and Jorge Urrutia}.
\newblock Crooked diagrams with few slopes.
\newblock \emph{Order}, 7(2):133--143, 1990{\natexlab{b}}.

\bibitem[{Dillencourt et~al.(2000)Dillencourt, Eppstein, and
  Hirschberg}]{DEH-JGAA00}
\textsc{Michael~B. Dillencourt, David Eppstein, and Daniel~S. Hirschberg}.
\newblock Geometric thickness of complete graphs.
\newblock \emph{J. Graph Algorithms Appl.}, 4(3):5--17, 2000.

\bibitem[{Ding and Oporowski(1995)}]{DO-JGT95}
\textsc{Guoli Ding and Bogdan Oporowski}.
\newblock Some results on tree decomposition of graphs.
\newblock \emph{J. Graph Theory}, 20(4):481--499, 1995.

\bibitem[{Ding and Oporowski(1996)}]{DO-DM96}
\textsc{Guoli Ding and Bogdan Oporowski}.
\newblock On tree-partitions of graphs.
\newblock \emph{Discrete Math.}, 149(1-3):45--58, 1996.

\bibitem[{Dujmovi{\'c} et~al.(2005{\natexlab{a}})Dujmovi{\'c}, Eppstein,
  Suderman, and Wood}]{DESW-PlanarSlopesSegments}
\textsc{Vida Dujmovi{\'c}, David Eppstein, Matthew Suderman, and David~R.
  Wood}.
\newblock Drawings of planar graphs with few slopes and segments.
\newblock Submitted, 2005{\natexlab{a}}.

\bibitem[{Dujmovi{\'c} et~al.(2005{\natexlab{b}})Dujmovi{\'c}, Morin, and
  Wood}]{DMW-SJC05}
\textsc{Vida Dujmovi{\'c}, Pat Morin, and David~R. Wood}.
\newblock Layout of graphs with bounded tree-width.
\newblock \emph{SIAM J. Comput.}, 34(3):553--579, 2005{\natexlab{b}}.

\bibitem[{Dujmovi{\'c} et~al.(2004)Dujmovi{\'c}, Suderman, and
  Wood}]{ReallyStraight-GD04}
\textsc{Vida Dujmovi{\'c}, Matthew Suderman, and David~R. Wood}.
\newblock Really straight graph drawings.
\newblock In \textsc{J\'{a}nos Pach}, ed., \emph{Proc. 12th International Symp.
  on Graph Drawing (GD '04)}, vol. 3383 of \emph{Lecture Notes in Comput.
  Sci.}, pp. 122--132. Springer, 2004.

\bibitem[{Dujmovi{\'c} and Wood(2004)}]{DujWoo-DMTCS04}
\textsc{Vida Dujmovi{\'c} and David~R. Wood}.
\newblock On linear layouts of graphs.
\newblock \emph{Discrete Math. Theor. Comput. Sci.}, 6(2):339--358, 2004.

\bibitem[{Dujmovi{\'c} and Wood(2006)}]{DujWoo-GD05}
\textsc{Vida Dujmovi{\'c} and David~R. Wood}.
\newblock Graph treewidth and geometric thickness parameters.
\newblock In \textsc{Patrick Healy and Nikola~S. Nikolov}, eds., \emph{Proc.
  13th International Symp. on Graph Drawing (GD '05)}, vol. 3843 of
  \emph{Lecture Notes in Comput. Sci.}, pp. 129--140. Springer, 2006.
\newblock \urlprefix\url{http://arxiv.org/math/0503553}.

\bibitem[{Duncan et~al.(2004)Duncan, Eppstein, and Kobourov}]{DEK-SoCG04}
\textsc{Christian~A. Duncan, David Eppstein, and Stephen~G. Kobourov}.
\newblock The geometric thickness of low degree graphs.
\newblock In \emph{Proc. 20th ACM Symp. on Computational Geometry (SoCG '04)},
  pp. 340--346. ACM Press, 2004.

\bibitem[{Eades et~al.(1997)Eades, Houle, and Webber}]{EHW97}
\textsc{Peter Eades, Michael~E. Houle, and Richard Webber}.
\newblock Finding the best viewpoints for three-dimensional graph drawings.
\newblock In \textsc{Giuseppe {Di~Battista}}, ed., \emph{Proc. 5th
  International Symp. on Graph Drawing (GD '97)}, vol. 1353 of \emph{Lecture
  Notes in Comput. Sci.}, pp. 87--98. Springer, 1997.

\bibitem[{Edenbrandt(1986)}]{Edenbrandt86}
\textsc{Anders Edenbrandt}.
\newblock Quotient tree partitioning of undirected graphs.
\newblock \emph{BIT}, 26(2):148--155, 1986.

\bibitem[{Eppstein(2004)}]{Eppstein-AMS}
\textsc{David Eppstein}.
\newblock Separating thickness from geometric thickness.
\newblock In \textsc{J\'{a}nos Pach}, ed., \emph{Towards a Theory of Geometric
  Graphs}, vol. 342 of \emph{Contemporary Mathematics}, pp. 75--86. Amer. Math.
  Soc., 2004.

\bibitem[{Ferber and J{\"u}rgensen(1970)}]{FJ70}
\textsc{K.~Ferber and Helmut J{\"u}rgensen}.
\newblock A programme for the drawing of lattices.
\newblock In \textsc{J.~Leech}, ed., \emph{Computational Problems in Abstract
  Algebra}, pp. 83--87. Pergamon, Oxford, 1970.

\bibitem[{Fomin and Golovach(2003)}]{FominGolovach-DAM03}
\textsc{Fedor~V. Fomin and Petr~A. Golovach}.
\newblock Interval degree and bandwidth of a graph.
\newblock \emph{Discrete Appl. Math.}, 129(2-3):345--359, 2003.

\bibitem[{Freese(2004)}]{Freese-04}
\textsc{Ralph Freese}.
\newblock Automated lattice drawing.
\newblock In \textsc{Peter~W. Eklund}, ed., \emph{Concept Lattices. Proc. 2nd
  International Conf. on Formal Concept Analysis (ICFCA '04)}, vol. 2961 of
  \emph{Lecture Notes in Comput. Sci.}, pp. 112--127. Springer, 2004.

\bibitem[{Halin(1991)}]{Halin91}
\textsc{Rudolf Halin}.
\newblock Tree-partitions of infinite graphs.
\newblock \emph{Discrete Math.}, 97:203--217, 1991.

\bibitem[{Houle and Webber(1998)}]{HW98}
\textsc{Michael~E. Houle and Richard Webber}.
\newblock Approximation algorithms for finding best viewpoints.
\newblock In \textsc{Sue Whitesides}, ed., \emph{Proc. 6th International Symp.
  on Graph Drawing (GD '98)}, vol. 1547 of \emph{Lecture Notes in Comput.
  Sci.}, pp. 210--223. Springer, 1998.

\bibitem[{Jamison(1984{\natexlab{a}})}]{Jamison84a}
\textsc{Robert~E. Jamison}.
\newblock Planar configurations which determine few slopes.
\newblock \emph{Geom. Dedicata}, 16(1):17--34, 1984{\natexlab{a}}.

\bibitem[{Jamison(1984{\natexlab{b}})}]{Jamison84}
\textsc{Robert~E. Jamison}.
\newblock Structure of slope-critical configurations.
\newblock \emph{Geom. Dedicata}, 16(3):249--277, 1984{\natexlab{b}}.

\bibitem[{Jamison(1986)}]{Jamison-DM86}
\textsc{Robert~E. Jamison}.
\newblock Few slopes without collinearity.
\newblock \emph{Discrete Math.}, 60:199--206, 1986.

\bibitem[{Jamison(1987)}]{Jamison87}
\textsc{Robert~E. Jamison}.
\newblock Direction trees.
\newblock \emph{Discrete Comput. Geom.}, 2(3):249--254, 1987.

\bibitem[{Kleitman and Pinchasi(2005)}]{KP-DCG05}
\textsc{Daniel~J. Kleitman and Rom Pinchasi}.
\newblock A note on caterpillar-embeddings with no two parallel edges.
\newblock \emph{Discrete Comput. Geom.}, 33(2):223--229, 2005.

\bibitem[{Malitz(1994)}]{Malitz94a}
\textsc{Seth~M. Malitz}.
\newblock Graphs with ${E}$ edges have pagenumber ${O}(\sqrt {E})$.
\newblock \emph{J. Algorithms}, 17(1):71--84, 1994.

\bibitem[{McKay(1985)}]{McKay-AC85}
\textsc{Brendan~D. McKay}.
\newblock Asymptotics for symmetric {$0$}-{$1$} matrices with prescribed row
  sums.
\newblock \emph{Ars Combin.}, 19(A):15--25, 1985.

\bibitem[{Onn and Pinchasi(2004)}]{OP-JCTA04}
\textsc{Shmuel Onn and Rom Pinchasi}.
\newblock A note on the minimum number of edge-directions of a convex polytope.
\newblock \emph{J. Combin. Theory Ser. A}, 107:147--151, 2004.

\bibitem[{Pach and P{\'a}lv{\"o}lgyi(2006)}]{PachPal-EJC06}
\textsc{J{\'a}nos Pach and D{\"o}m{\"o}t{\"o}r P{\'a}lv{\"o}lgyi}.
\newblock Bounded-degree graphs can have arbitrarily large slope numbers.
\newblock \emph{Electron. J. Combin.}, 13(1):N1, 2006.

\bibitem[{Pach et~al.(2004{\natexlab{a}})Pach, Pinchasi, and
  Sharir}]{PPS-JCTA04}
\textsc{J{\'a}nos Pach, Rom Pinchasi, and Micha Sharir}.
\newblock On the number of directions determined by a three-dimensional points
  set.
\newblock \emph{J. Combin. Theory Ser. A}, 108(1):1--16, 2004{\natexlab{a}}.

\bibitem[{Pach et~al.(2004{\natexlab{b}})Pach, Pinchasi, and
  Sharir}]{PPS-SoCG04}
\textsc{J\'{a}nos Pach, Rom Pinchasi, and Micha Sharir}.
\newblock Solution of {S}cott's problem on the number of directions determined
  by a point set in 3-space.
\newblock In \emph{Proc. 20th ACM Symp. on Computational Geometry (SoCG '04)},
  pp. 76--85. ACM Press, 2004{\natexlab{b}}.

\bibitem[{Reed(2003)}]{Reed-AlgoTreeWidth03}
\textsc{Bruce~A. Reed}.
\newblock Algorithmic aspects of tree width.
\newblock In \textsc{Bruce~A. Reed and Cl{\'a}udia~L. Sales}, eds.,
  \emph{Recent Advances in Algorithms and Combinatorics}, pp. 85--107.
  Springer, 2003.

\bibitem[{Scheffler(1989)}]{Scheffler89}
\textsc{Petra Scheffler}.
\newblock \emph{Die Baumweite von Graphen als ein Ma{\ss} f{\"u}r die
  Kompliziertheit algorithmischer Probleme}.
\newblock Ph.D. thesis, Akademie der Wissenschaften der DDR, Berlin, Germany,
  1989.

\bibitem[{Scott(1970)}]{Scott70}
\textsc{Paul~R. Scott}.
\newblock On the sets of directions determined by {$n$} points.
\newblock \emph{Amer. Math. Monthly}, 77:502--505, 1970.

\bibitem[{Seese(1985)}]{Seese85}
\textsc{Detlef Seese}.
\newblock Tree-partite graphs and the complexity of algorithms.
\newblock In \textsc{Lothar Budach}, ed., \emph{Proc. International Conf. on
  Fundamentals of Computation Theory}, vol. 199 of \emph{Lecture Notes in
  Comput. Sci.}, pp. 412--421. Springer, 1985.

\bibitem[{Suderman(2004)}]{Suderman-IJCGA04}
\textsc{Matthew Suderman}.
\newblock Pathwidth and layered drawings of trees.
\newblock \emph{Internat. J. Comput. Geom. Appl.}, 14(3):203--225, 2004.

\bibitem[{Ungar(1982)}]{Ungar-JCTA82}
\textsc{Peter Ungar}.
\newblock {$2N$} noncollinear points determine at least {$2N$} directions.
\newblock \emph{J. Combin. Theory Ser. A}, 33(3):343--347, 1982.

\bibitem[{Wade and Chu(1994)}]{WadeChu-CJ94}
\textsc{Greg~A. Wade and Jiang-Hsing Chu}.
\newblock Drawability of complete graphs using a minimal slope set.
\newblock \emph{The Computer Journal}, 37(2):139--142, 1994.

\bibitem[{Wood(2006{\natexlab{a}})}]{Wood-IntersectionGraphOrderings}
\textsc{David~R. Wood}.
\newblock Characterisations of intersection graphs by vertex orderings.
\newblock \emph{Australas. J. Combin.}, 34:261--268, 2006{\natexlab{a}}.

\bibitem[{Wood(2006{\natexlab{b}})}]{Wood-TreePartitions}
\textsc{David~R. Wood}.
\newblock A note on tree-partition-width, 2006{\natexlab{b}}.
\newblock \urlprefix\url{http://arxiv.org/math/0602507}.

\bibitem[{Wormald(1978)}]{Wormald78}
\textsc{Nicholas Wormald}.
\newblock \emph{Some problems in the enumeration of labelled graphs}.
\newblock Ph.D. thesis, Newcastle-upon-Tyne, United Kingdom, 1978.

\end{thebibliography}

\def\soft#1{\leavevmode\setbox0=\hbox{h}\dimen7=\ht0\advance \dimen7
  by-1ex\relax\if t#1\relax\rlap{\raise.6\dimen7
  \hbox{\kern.3ex\char'47}}#1\relax\else\if T#1\relax
  \rlap{\raise.5\dimen7\hbox{\kern1.3ex\char'47}}#1\relax \else\if
  d#1\relax\rlap{\raise.5\dimen7\hbox{\kern.9ex \char'47}}#1\relax\else\if
  D#1\relax\rlap{\raise.5\dimen7 \hbox{\kern1.4ex\char'47}}#1\relax\else\if
  l#1\relax \rlap{\raise.5\dimen7\hbox{\kern.4ex\char'47}}#1\relax \else\if
  L#1\relax\rlap{\raise.5\dimen7\hbox{\kern.7ex
  \char'47}}#1\relax\else\message{accent \string\soft \space #1 not
  defined!}#1\relax\fi\fi\fi\fi\fi\fi} \def\cprime{$'$}

\end{document}